\documentclass[12pt]{article}

\usepackage{arxiv}

\usepackage{amsmath, amssymb, amsfonts}
\usepackage{amsthm}

\newtheoremstyle{estiloTeorema}
  {12pt}                
  {12pt}                
  {\itshape}            
  {}                    
  {\bfseries}           
  {.}                   
  {.5em}                
  {}                    

\newtheoremstyle{estiloDefinicion}
  {12pt}                
  {12pt}                
  {\normalfont}         
  {}                    
  {\bfseries}           
  {.}                   
  {.5em}                
  {}                    

\theoremstyle{estiloTeorema}
\newtheorem{theorem}{Theorem}[section]

\theoremstyle{estiloDefinicion} 
\newtheorem{definition}{Definition}[section]

\usepackage[utf8]{inputenc} 
\usepackage[T1]{fontenc}    
\usepackage{hyperref}       
\usepackage{url}            
\usepackage{booktabs}       
\usepackage{amsfonts}       
\usepackage{nicefrac}       
\usepackage{microtype}      
\usepackage{lipsum}		
\usepackage{graphicx}
\usepackage{natbib}
\usepackage{doi}
\usepackage{amsmath, amssymb, amsfonts}
\usepackage{amsthm}

\title{The Schwarz-Pick Lemma as a Consequence of the Ahlfors-Schwarz-Pick Lemma}

\author{ {\hspace{1mm}Rafael Benjumea Cejas} \\
	Departamento de Análisis Matemático\\
	Universidad de Sevilla\\
	Seville, Spain \\
	\texttt{rafbencej@alum.us.es} \\
	\And
	{\hspace{1mm}Juan Carlos García Vázquez} \\
	Departamento de Análisis Matemático\\
	Universidad de Sevilla\\
	Seville, Spain\\
	\texttt{garcia@us.es} \\
}

\date{}

\renewcommand{\headeright}{}
\renewcommand{\undertitle}{}
\renewcommand{\shorttitle}{The Schwarz-Pick Lemma as a Consequence of the Ahlfors-Schwarz-Pick Lemma}

\hypersetup{
pdftitle={The Schwarz-Pick Lemma as a Consequence of the Ahlfors-Schwarz-Pick Lemma},
pdfsubject={q-bio.NC, q-bio.QM},
pdfauthor={Rafael Benjumea Cejas, Juan Carlos García Vázquez},
pdfkeywords={Schwarz-Pick Lemma, Ahlfors-Schwarz-Pick Lemma, Hyperbolic Metric, Conformal Metrics},
}

\begin{document}
\maketitle

\begin{abstract}
	The aim of this article is to give an elementary proof of the fact that the Schwarz-Pick Lemma follows from the Ahlfors-Schwarz-Pick Lemma
\end{abstract}


\section{Introduction}
In 1907, Carathéodory introduced the version of the Schwarz Lemma that it is known today. Let us recall it.
\begin{theorem}[Schwarz Lemma]
    Let $f:\mathbb{D} \longrightarrow \mathbb{D}$ be a holomorphic function with $f(0)=0$. The following inequalities hold:
    \begin{itemize}
        \item $|f(z)| \leq |z|$ for all $z \in \mathbb{D}$.
        \item $|f'(0)| \leq 1$.
    \end{itemize}
    Moreover, the following are equivalent:
    \begin{enumerate}
        \item There exists $z_0 \in \mathbb{D}$ such that $|f(z_0)|=|z_0|$ or $|f'(0)|=1$.
        \item $f$ is a rotation, that is to say,
        \begin{equation*}
            f(z)=cz \ \ \forall z \in \mathbb{D},
        \end{equation*}
        where $c \in \partial \mathbb{D}$.
        \item $|f(z)|=|z|$ for all $z \in \mathbb{D}$ or $|f'(0)|=1$.
    \end{enumerate}
\end{theorem}

The Schwarz Lemma is a result about the distorsion of holomorphic self-maps of the unit disk. Despite its simplicity, the Schwarz Lemma has far-reaching generalizations. The most celebrated one is the Schwarz-Pick Lemma, established by Pick in 1916.

\begin{theorem}[Schwarz-Pick Lemma]
    Let $f:\mathbb{D} \longrightarrow \mathbb{D}$ be a holomorphic function. For any $z,w \in \mathbb{D}$, the following inequalities hold:
    \begin{equation*}
        d_\mathbb{D}(f(z),f(w)) \leq d_\mathbb{D}(z,w), \qquad \lambda_\mathbb{D}(f(z))|f'(z)| \leq \lambda_\mathbb{D}(z).
    \end{equation*}
    Moreover, the following are equivalent:
    \begin{enumerate}
        \item There exist $z_0,w_0 \in \mathbb{D}$ such that
        \begin{equation*}
            d_\mathbb{D}(f(z_0),f(w_0))=d_\mathbb{D}(z_0,w_0), \qquad \lambda_\mathbb{D}(f(z_0))|f'(z_0)| = \lambda_\mathbb{D}(z_0).
        \end{equation*}
        \item $f$ is a biholomorphic map (i.e, $f$ is one-to-one, onto and holomorphic).
        \item For all $z,w \in \mathbb{D}$,
        \begin{equation*}
            d_\mathbb{D}(f(z),f(w))=d_\mathbb{D}(z,w), \qquad \lambda_\mathbb{D}(f(z))|f'(z)|=\lambda_\mathbb{D}(z).
        \end{equation*}
    \end{enumerate}
\end{theorem}

In the previous result, $d_\mathbb{D}$ denotes the Poincaré distance (Theorem 3.2), whereas $\lambda_\mathbb{D}$ is the hyperbolic metric on $\mathbb{D}$ (Definition 2.3). Taking a holomorphic function $f:\mathbb{D} \longrightarrow \mathbb{D}$ with $f(0)=0$ and making $z=0$ in the inequalities of the Schwarz-Pick Lemma, it follows from the formula of $d_\mathbb{D}$ and the definition of $\lambda_\mathbb{D}$ that the Schwarz Lemma is a consequence of the Schwarz-Pick Lemma. 

The Schwarz-Pick lemma states that, with respect to the Poincaré distance, every holomorphic self-map of the unit disk is distance-decreasing. Note that the hypothesis $f(0)=0$ is cast away. Also, observe that the focus of the Schwarz-Pick Lemma is not the holomorphic map $f$ but the hyperbolic metric $\lambda_\mathbb{D}$. Thus, the nature of the Schwarz Lemma changes from analytical to geometrical.

It turns out that the language in whch the Schwarz Lemma reaches all its generality is that of conformal metrics (Definition 2.1) and their curvature (Definition 2.8). It was Ahlfors who first adopted this pioneering perspective, providing the following result in 1938.

\begin{theorem}[Maximality principle of the hyperbolic metric]
    Let $\rho$ be a conformal metric on $\mathbb{D}$ such that $K_\rho(z) \leq -1$ for all $z \in \mathbb{D}$. Then
    \begin{equation*}
        \rho(z) \leq \lambda_\mathbb{D}(z) \ \ \forall z \in \mathbb{D}.
    \end{equation*}
\end{theorem}

For a proof of this result, consult \cite{Beardon-Minda}. Thus, out of all conformal metrics on $\mathbb{D}$ with curvature at most $-1$, it is the hyperbolic metric that has greatest value at each point. Later, in 1962, Heins improved Ahlfors' result by characterizing the equality case. As a consequence, he obtained the Ahlfors-Schwarz-Pick Lemma.

\begin{theorem}[Ahlfors-Schwarz-Pick Lemma]
    Let $\rho$ be a conformal metric on $\mathbb{D}$ such that $K_\rho(z) \leq -1$ for all $z \in \mathbb{D}$. The following inequality holds:
    \begin{equation*}
        \rho(z) \leq \lambda_\mathbb{D}(z) \ \ \forall z \in \mathbb{D}.
    \end{equation*}
    Moreover, the following are equivalent:
    \begin{enumerate}
        \item There exists $z_0 \in \mathbb{D}$ such that $\rho(z_0)=\lambda_\mathbb{D}(z_0)$.
        \item For all $z \in \mathbb{D}$, we have that $\rho(z)=\lambda_\mathbb{D}(z)$. 
    \end{enumerate}
\end{theorem}

An elementary proof of the equality case may be found in \cite{Kraus-Roth}. Unlike the Schwarz Lemma or the Schwarz-Pick Lemma, there is not any (holomorphic) function in the hypothesis. Therefore, the Ahlfors-Schwarz-Pick Lemma is purely geometrical, namely, differential-geometrical.

It seems that there is no relationship between the Schwarz-Pick Lemma and the Ahlfors-Schwarz-Pick Lemma other than their names. In \cite{Beardon-Minda}, A. F. Beardon and C. D. Minda affirm that the Schwarz-Pick Lemma is a consequence of the Ahlfors-Schwarz-Pick Lemma, but no explicit proof is given. The aim of this paper is to provide an elementary demonstration of this fact. Namely, we will prove the following theorem.

\begin{theorem}
    The Schwarz-Pick Lemma is a consequence of the Ahlfors-Schwarz-Pick Lemma.
\end{theorem}

\section{Conformal metrics}
\label{sec:headings}
In this section, we introduce the differential-geometric concepts concerning the Ahlfors-Schwarz-Pick Lemma, which will be needed in the proof of Theorem 1.5.
\begin{definition}
    Let $\Omega \subseteq \mathbb{C}$ be a domain. A \textit{conformal metric} on $\Omega$ is a function 
    \begin{equation*}
        \rho: \Omega \longrightarrow [0,+\infty) 
    \end{equation*}
    which verifies the following:
    \begin{itemize}
        \item $\rho$ is continuous.
        \item $\rho$ vanishes on a discrete subset of $\Omega$.
        \item $\rho$ is of class $\mathcal{C}^2$ at every point $z \in \Omega$ such that $\rho(z)>0$.
    \end{itemize}
\end{definition}

\begin{definition}
    Let $\Omega \subseteq \mathbb{C}$ be a domain and $\rho$ a conformal metric on $\Omega$. The \textit{length of $\xi \in \mathbb{C}$ at $z \in \Omega$ with respect to $\rho$} is
    \begin{equation*}
        \|\xi\|_{\rho,z}=\rho(z)|\xi|.
    \end{equation*}
\end{definition}

Essentially, a conformal metric on a domain $\Omega \subseteq \mathbb{C}$ is a function that provides a new way to measure the length of a vector $\xi \in T_z\Omega \cong \mathbb{C}$ from its euclidean norm. Here, $T_z\Omega$ denotes the tangent plane of $\Omega$ at $z$. One of the most important examples of conformal metrics is the hyperbolic metric on the unit disk $\mathbb{D}$.

\begin{definition}
    The \textit{hyperbolic metric on the unit disk $\mathbb{D}$} or the \textit{Poincaré metric} is
    \begin{equation*}
        \lambda_\mathbb{D}(z)=\dfrac{2}{1-|z|^2} \ \ \forall z \in \mathbb{D}.
    \end{equation*}
\end{definition}

It is not difficult to check that the hyperbolic metric is indeed a conformal metric on $\mathbb{D}$. Conformal metrics allow us to measure the length of piecewise regular curves.

\begin{definition}
    Let $\Omega \subseteq \mathbb{C}$ be a domain, $\rho$ a conformal metric on $\Omega$ and $\gamma:[a,b] \longrightarrow \Omega$ a piecewise regular curve. The \textit{$\rho-$length} of $\gamma$ is
    \begin{equation*}
        \ell_\rho(\gamma)=\int_a^b \|\gamma'(t)\|_{\rho,\gamma(t)}\,dt = \int_a^b \rho(\gamma(t))|\gamma'(t)|\, dt.
    \end{equation*}
\end{definition}

The preceding definition does not depend on the parametrization of the curve $\gamma$. Besides, it enables us to introduce a distance function on $\Omega$.

\begin{definition}
    Let $\Omega \subseteq \mathbb{C}$ be a domain, $\rho$ a conformal metric on $\Omega$ and $z,w \in \Omega$. The \textit{$\rho-$distance between $z$ and $w$} is
    \begin{equation*}
        d_\rho(z,w)=\inf\{\ell_{\rho}(\eta) : \eta \in \mathcal{C}_\Omega(z,w)\},
    \end{equation*}
    where $\mathcal{C}_\Omega(z,w)$ is the set of all piecewise regular curves in $\Omega$ joining $z$ and $w$.
\end{definition}

For a proof of the fact that $d_\rho$ defines a distance on $\Omega$, consult \cite{Carroll}. The maps that preserve the distance functions associated to conformal metrics are relevant. 

\begin{definition}
    Let $\Omega_1, \Omega_2 \subseteq \mathbb{C}$ be domains, $\rho_1, \rho_2$ conformal metrics on $\Omega_1$ and $\Omega_2$, respectively, and $f:\Omega_1 \longrightarrow \Omega_2$ an onto function. We say that $f$ is an \textit{isometry with respect to $\rho_1$ and $\rho_2$} if
    \begin{equation*}
        d_{\rho_1}(z,w)=d_{\rho_2}(f(z),f(w)) \ \ \forall z,w \in \Omega_1.
    \end{equation*}
\end{definition}

The following definition will be helpful in the study of isometries.

\begin{definition}
    Let $\Omega_1, \Omega_2 \subseteq \mathbb{C}$ be domains, $\rho$ a conformal metric on $\Omega_2$ and $f:\Omega_1 \longrightarrow \Omega_2$ a holomorphic function. The \textit{pullback of $\rho$ by $f$} is 
    \begin{equation*}
        f^*\rho(z)=\|f'(z)\|_{\rho,f(z)}=\rho(f(z))|f'(z)| \ \ \forall z \in \Omega_1.
    \end{equation*}
\end{definition}

If the function $f$ in the previous definition is non-constant, then the pullback of $\rho$ by $f$ is a conformal metric on $\Omega_1$. A proof of this fact can be found in \cite{Krantz}. Thanks to the pullback of conformal metrics, we can provide a characterization of isometries.

\begin{theorem}[Characterization of isometries]
    Let $\Omega_1, \Omega_2 \subseteq \mathbb{C}$ be domains, $\rho_1, \rho_2$ conformal metrics on $\Omega_1$ and $\Omega_2$, respectively, and $f:\Omega_1 \longrightarrow \Omega_2$ a biholomorphic map. The following are equivalent:
    \begin{enumerate}
    \item $f$ is an isometry with respect to $\rho_1$ and $\rho_2$.
    \item The pullback of $\rho_2$ by $f$ is precisely $\rho_1$, that is to say, the following equality holds:
    \begin{equation*}
        f^*\rho_2(z)=\rho_2(f(z))|f'(z)|=\rho_1(z) \ \ \forall z \in \Omega_1.
    \end{equation*}
    \item The length of piecewise regular curves is preserved by $f$, namely,
    \begin{equation*}
        \ell_{\rho_1}(\gamma)=\ell_{\rho_2}(f \circ \gamma)
    \end{equation*}
    for all piecewise regular curves $\gamma$ in $\Omega_1$.
    \end{enumerate}
\end{theorem}

A proof of the previous characterization may also be found in \cite{Krantz}. Finally, we review the notion of curvature of a conformal metric.

\begin{definition}
    Let $\Omega \subseteq \mathbb{C}$ be a domain, $\rho$ a conformal metric on $\Omega$ and $z \in \Omega$ with $\rho(z)>0$. The \textit{curvature of $\rho$ at $z$} is
    \begin{equation*}
        K_\rho(z)=-\dfrac{\Delta \log \rho(z)}{\rho(z)^2},
    \end{equation*}
    where $\Delta=\partial_{xx}^{2}+\partial_{yy}^{2}$ denotes the laplacian.
\end{definition}

Under suitable conditions, the pullback of a conformal metric does not change its curvature.

\begin{theorem}[Invariance of curvature]
    Let $\Omega_1,\Omega_2 \subseteq \mathbb{C}$ be domains, $\rho$ a conformal metric on $\Omega_2$ and $f:\Omega_1 \longrightarrow \Omega_2$ a nonconstant holomorphic function. The equality
    \begin{equation*}
        K_{\Omega_1,f^*\rho}(z)=K_{\Omega_2,\rho}(f(z))
    \end{equation*}
    holds for all $z \in \Omega_1$ such that $f'(z) \neq 0$ and $\rho(f(z)) >0$.
\end{theorem}

A proof of this theorem may be found in \cite{Beardon-Minda}.

\section{The hyperbolic metric on $\mathbb{D}$}
\label{sec:others}

 In this section, we state some facts about the Poincaré metric that will be necessary to prove Theorem 1.5. We denote by $\textnormal{Aut}(\mathbb{D})$ the set of all automorphisms of $\mathbb{D}$, that is to say, biholomorphic self-maps of $\mathbb{D}$. We recall that it consists of the functions of the form
\begin{equation*}
   f(z)=e^{i\theta}\dfrac{z-a}{1-\overline{a}z} \ \ \forall z \in \mathbb{D}, 
\end{equation*}
where $\theta \in \mathbb{R}$ and $a \in \partial \mathbb{D}$.

\begin{theorem}
    Every automorphism of $\mathbb{D}$ is an isometry with respect to the Poincaré metric.
\end{theorem}

For a proof of this result, consult \cite{Krantz}. Thanks to the previous theorem, we are able to give an explicit formula for the distance function asociated to the hyperbolic metric on $\mathbb{D}$. We write $d_\mathbb{D}=d_{\lambda_\mathbb{D}}$ to denote the $\lambda_\mathbb{D}-$distance on $\mathbb{D}$ and we call it \textit{hyperbolic distance on $\mathbb{D}$} or \textit{Poincaré distance}. 

\begin{theorem}
   The hyperbolic distance on $\mathbb{D}$ obeys the following formula:
\begin{equation*}
    d_\mathbb{D}(z,w)=\log \dfrac{1+\dfrac{z-w}{1-\overline{w}z}}{1-\dfrac{z-w}{1-\overline{w}z}} \ \ \forall z, w \in \mathbb{D}.
\end{equation*} 
\end{theorem}

The reference \cite{Beardon-Minda} contains a proof of the previous result.

\section{Proof of Theorem 1.5}
Finally, we prove Theorem 1.5. Let $f:\mathbb{D} \longrightarrow \mathbb{D}$ be a nonconstant holomorphic function. According to Theorem 2.2, the conformal metric $f^*\lambda_\mathbb{D}$ has curvature identically $-1$. Thus, by the Ahlfors-Schwarz-Pick Lemma,
we have 
\begin{equation*}
    f^*\lambda_\mathbb{D}(z) \leq \lambda_\mathbb{D}(z) \ \ \forall z \in \mathbb{D}.
\end{equation*}
By Definition 2.7, this inequality is the same as
\begin{equation*}
    \lambda_\mathbb{D}(f(z))|f'(z)| \leq \lambda_\mathbb{D}(z) \ \ \forall z \in \mathbb{D}.
\end{equation*}
This is precisely one of the inequalities in the Schwarz-Pick Lemma. The other one is easily deduced by integration over arbitrary piecewise regular curves in $\mathbb{D}$.

Now, let us study the equality case. Once again by Definition 12, there exists $z_0 \in \mathbb{D}$ with \[ \lambda_\mathbb{D}(f(z_0))|f'(z_0)|=\lambda_\mathbb{D}(z_0)\] if and only if \[f^*\lambda_\mathbb{D}(z_0)=\lambda_\mathbb{D}(z_0).\] By the Ahlfors-Schwarz-Pick Lemma, the previous equality holds if and only if \[f^*\lambda_\mathbb{D}(z)=\lambda_\mathbb{D}(z) \ \ \forall z \in \mathbb{D}.\] and, in turn, if and only if \[\lambda_\mathbb{D}(f(z))|f'(z)|=\lambda_\mathbb{D}(z) \ \ \forall z \in \mathbb{D}.\] Thus, it remains to prove that the last equality holds if and only if $f \in \textnormal{Aut}(\mathbb{D})$.

If $f \in \textnormal{Aut}(\mathbb{D})$, then Theorem 3.1 guarantees that $f$ is an isometry with respect to $\lambda_\mathbb{D}$. Therefore, Theorem 2.1 implies that
\begin{equation*}
    \lambda_\mathbb{D}(f(z))|f'(z)|=\lambda_\mathbb{D}(z) \ \ \forall z \in \mathbb{D}.
\end{equation*}

Next, let us prove the converse. First, we demonstrate that $f$ is onto. Suppose that $f(\mathbb{D}) \subsetneq \mathbb{D}$. In particular, $\overline{f(\mathbb{D})} \subsetneq \overline{\mathbb{D}}$, from where we deduce that $\overline{\mathbb{D}} \setminus \overline{f(\mathbb{D})}$ is not the empty set. Now, we take $z_0 \in f(\mathbb{D})$ and $z_1 \in \overline{\mathbb{D}} \setminus \overline{f(\mathbb{D})}$ and consider a continuous curve $\eta:[0,1] \longrightarrow \overline{\mathbb{D}}$ with $\eta(0)=z_0$ and $\eta(1)=z_1$. There exist $0<r<s<1$ such that $\eta([0,r]) \subseteq f(\mathbb{D})$ and $\eta([s,1]) \subseteq \overline{\mathbb{D}} \setminus \overline{f(\mathbb{D})}$. It is easy to check that
\begin{equation*}
    \overline{\mathbb{D}}=f(\mathbb{D}) \sqcup \partial f(\mathbb{D}) \sqcup (\overline{\mathbb{D}} \setminus \overline{f(\mathbb{D})}). 
\end{equation*}
Let us denote $\eta^*=\eta([0,1])$. If $\eta^* \cap \partial f(\mathbb{D}) = \emptyset$, then $\eta^* \subseteq f(\mathbb{D}) \sqcup \overline{\mathbb{D}} \setminus \overline{f(\mathbb{D})}$, which would imply that $\eta^*$ is not connected. However, this is absurd since $\eta^*$ is connected as $\eta$ is continuous and $[0,1]$ is connected. Thus, $\partial f(\mathbb{D}) \cap \eta^* \neq \emptyset$.

Given $w \in \partial f(\mathbb{D}) \cap \mathbb{D}$ and a sequence $(z_n)_{n=1}^{\infty}$ verifying $g(z_n) \longrightarrow w$, it follows that $|z_n| \longrightarrow 1$. Indeed, otherwise, there would exist $R \in (0,1)$ such that $|z_n|<R$ for all $n \in \mathbb{N}$. Since the disk $\overline{D}(0,R)$ is compact, the Bolzano-Weierstrass Theorem assures the existence of $z \in \overline{D}(0,R)$ and a subsequence $(z_{n_k})_{k=1}^{\infty}$ with $z_{n_k} \longrightarrow z$. As a consequence, $f(z_{n_k}) \longrightarrow f(z)$, so $f(z)=w$. By the Open Mapping Theorem, $f(\mathbb{D})$ is open in $\mathbb{D}$, so $f(z) \in \mathbb{D}$. On the other hand, $w \in \partial f (\mathbb{D})$, from where we get a contradiction. Thus, $|z_n| \longrightarrow 1$, as we desired to prove.

Now, we observe that
\begin{equation*}
    \lambda_\mathbb{D}(f(z_n))|f'(z_n)|=\dfrac{2|f'(z_n)|}{1-|f(z_n)|^2}=\dfrac{2}{1-|z_n|^2}=\lambda_\mathbb{D}(z_n) \ \ \forall n \in \mathbb{N}.
\end{equation*}
Due to the fact that $|z_n| \longrightarrow 1$, we have
\begin{equation*}
    \underset{n \rightarrow \infty}{\textnormal{lim}} \dfrac{2}{1-|z_n|^2}=+\infty.
\end{equation*}
Furthermore, from Cauchy's estimates we get the following upper bound:
\begin{equation*}
    |f'(z_n)| \leq \dfrac{1}{r} \cdot \underset{|s-z_n|=r}{\textnormal{max}} |f(s)| \leq \dfrac{1}{r} \ \ \forall n \in \mathbb{N}.
\end{equation*}
Here, $r>0$ is such that $\overline{D}(z_n,r) \subseteq \mathbb{D}$ for all $n \in \mathbb{N}$ (for instance, $r=d(\partial f(\mathbb{D}),\partial \mathbb{D})/2$). Consequently,
\begin{equation*}
    \underset{n \rightarrow \infty}{\textnormal{lim}} \dfrac{2|f'(z_n)|}{1-|f(z_n)|^2} \leq \dfrac{1}{r} \underset{n \rightarrow \infty}{\textnormal{lim}} \dfrac{2}{1-|f(z_n)|^2} = \dfrac{2}{r(1-|w|^2)} < +\infty.
\end{equation*}
This is a contradiction. In conclusion, $f(\mathbb{D})=\mathbb{D}$, that is to say, $f$ is onto.

Lastly, we prove that $f$ is one-to-one. Since $\lambda_\mathbb{D}(z)>0$ for all $z \in \mathbb{D}$, it follows that $f'(z) \neq 0$ for all $z \in \mathbb{D}$ and, therefore, $f$ is locally univalent. Thus, by the Inverse Function Theorem, given $z_0 \in \mathbb{D}$ and $w_0:=f(z_0)\in f(\mathbb{D})=\mathbb{D}$, there exist open sets $U_0,V_0 \subseteq \mathbb{D}$ such that $z_0 \in U_0$, $w_0 \in V_0$ and $f:U_0 \longrightarrow V_0$ is a biholomorphic map. If $g_0:V_0 \longrightarrow U_0$ is the local inverse of $f$ in $V_0$, we see that $f \circ g_0 = \textnormal{id}_{V_0}$ and $g_0 \circ f = \textnormal{id}_{U_0}$. In these equalities, we understand that $f$ is restricted to $U_0$.

Now, let $\gamma:[0,1] \longrightarrow \mathbb{D}$ be a continuous curve with $\gamma(0)=w_0$. Given $t \in [0,1]$, if $W_t$ is an open neighborhood of $\gamma(t)$, then $\{W_t : t \in [0,1]\}$ is an open covering of $\gamma^*$. Since $\gamma^*$ is compact, there exist $0=t_0<t_1< \cdots < t_n=1$ verifying
\begin{equation*}
    \gamma^* \subseteq W_0 \cup W_1 \cup \cdots \cup W_n,
\end{equation*}
where $W_k:=W_{t_k}$ for all $k=0,1,...,n$. Without loss of generality, we may suppose that $W_0=V_0$. We affirm that $W_k \cap W_{k+1} \neq \emptyset$ for all $k=0,1,...,n$. In effect, otherwise, $\gamma^*$ would not be a connected set, which is absurd.

Next, by the Inverse Function Theorem, there exist open sets $U_1,...,U_n \subseteq f(\mathbb{D}) = \mathbb{D}$ such that the functions
\begin{equation*}
    g_1:W_1 \longrightarrow U_1, \quad g_2:W_2 \longrightarrow U_2, \quad ..., \quad g_n:W_n \longrightarrow U_n
\end{equation*}
are local inverses of $f$. In $V_0 \cap W_1$, both $g_0$ and $g_1$ are local inverses of $f$, so
\begin{equation*}
    f(g_0(w))=w=f(g_1(w)) \ \ \forall w \in V_0 \cap W_1.
\end{equation*}
Therefore, since $f$ is locally univalent, we deduce $g_0(w)=g_1(w)$ for all $w \in V_0 \cap W_1$. Moreover, by analytic continuation, both $g_0$ and $g_1$ extend to a function $G_{0,1}$, which is holomorphic in $V_0 \cup W_1$. The function $G_{0,1}$ is a local inverse of $f$ in $V_0 \cup W_1$. In effect, given $v \in V_0$ arbitrary but fixed,
\begin{equation*}
    (G_{0,1} \circ f)(v)=(g_0 \circ f)(v) = v \quad \textnormal{y} \quad (f \circ G_{0,1})(v)=(f \circ g_0)(v)=v.
\end{equation*}
The same holds if $v \in W_1$ or $v \in V_0 \cap W_1$. Thanks to a similar reasoning, we have $G_{0,1}(w)=g_2(w)$ for all $w \in W_1 \cap W_2$. Thus, once again by analytic continuation, $G_{0,1}$ and $g_2$ extend to a function $G_{0,1,2}$, holomorphic in $V_0 \cup W_1 \cup W_2$, which is a local inverse of $g$ in $V_0 \cup W_1 \cup W_2$. Continuing this process, we elucidate that $g_0$ analitically extends to a function $G$ holomorphic in the tubular neigborhood of $\gamma^*$ given by
\begin{equation*}
    B:=V_0 \cup W_1 \cup W_2 \cup \cdots \cup W_n.
\end{equation*}
Besides, $G$ is a local inverse of $f$ in $B$. Thus, since $\mathbb{D}$ is a simply connected domain and $\gamma$ is arbitrary, the Monodromy Theorem guarantees that $g_0$ extends to a holomorphic function in $\mathbb{D}$, which we denote again by $G$.

Finally, let us see that $G$ is the global inverse of $f$. Let $z \in \mathbb{D}$ arbitrary but fixed and $\alpha:[0,1] \longrightarrow \mathbb{D}$ a continuous curve joining $w_0$ and $z$. Following the same steps as above, there exists a tubular neighborhood $A \subseteq \mathbb{D}$ of $\alpha^*$ and a function $k$ which is holomorphic in $A$ and an analytic continuation of $g_0$, as well as a local inverse of $f$ in $A$. Therefore,
\begin{equation*}
    (G \circ f)(z)=(k \circ f)(z)=z \quad \textnormal{y} \quad (f \circ G)(z)=(f \circ k)(z)=z,
\end{equation*}
so $G$ is the global inverse of $f$. In conclusion, $f \in \textnormal{Aut}(\mathbb{D})$, as we desired to prove.

\bibliographystyle{unsrtnat}
\bibliography{references}  






\end{document}